\documentstyle[12pt]{article}
\parskip=\medskipamount                 
\thicklines                     
\newcommand{\bimn}[7]{\bibitem{#1}#2,
{\em #3},
{ #4}$\;${\bf
#5}$\;$(#6)$\;${#7}.}

\def\inbar{\vrule height1.5ex width.4pt depth0pt}
\def\IN{\relax{\rm I\kern-.18em N}}
\def\IQ{\relax\,\hbox{$\inbar\kern-.3em{\rm Q}$}}
\def\IR{\relax{\rm I\kern-.18em R}}
\def\ZZ{\relax{\sf Z\kern-.4em Z}}
\def\a{\alpha} \def\b{\beta}   \def\e{\epsilon} 
 \def\l{\lambda} 
  \def\s{\sigma}

 \def\cK{{\cal K}} \def\cL{{\cal L}} 
 \def\cO{{\cal O}}  
 
\newtheorem{theorem}{Theorem}[section]
\newtheorem{proposition}{Proposition}[section]

\newtheorem{conjecture}{Conjecture}[section]
\newtheorem{lemma}{Lemma}[section]

\catcode`\@=11

\newif\if@fewtab\@fewtabtrue

\catcode`\@=11

\newif\if@fewtab\@fewtabtrue

{\count255=\time\divide\count255 by 60
\xdef\hourmin{\number\count255}
\multiply\count255 by-60\advance\count255 by\time
\xdef\hourmin{\hourmin:\ifnum\count255<10 0\fi\the\count255}}
\def\ps@draft{\let\@mkboth\@gobbletwo
    \def\@oddhead{}
    \def\@oddfoot
      {\hbox to 7 cm{\footnotesize {\em Draft of \jobname:} \draftdate
       \hfil}\hskip -7cm\hfil\rm\thepage \hfil}
    \def\@evenhead{}\let\@evenfoot\@oddfoot}

\def\ceqno{\global\@fewtabfalse
    \ifcase\@eqcnt \def\@tempa{& & &}\or \def\@tempa{& &}
      \or \def\@tempa{&}
      \or\def\@tempa{}\fi\@tempa
{\rm(\theequation)}}

\def\aeqno#1{\global\@fewtabfalse
    \ifcase\@eqcnt \def\@tempa{& & &}\or \def\@tempa{& &}
      \or \def\@tempa{&}
      \or\def\@tempa{}\fi\@tempa
{\rm(\theequation,#1)}}

\def\label#1{\ifnum\draftcontrol=1
 \global\def\draftnote{$\scriptstyle #1$}\fi
 \@bsphack\if@filesw {\let\thepage\relax
   \def\protect{\noexpand\noexpand\noexpand}%
\xdef\@gtempa{\write\@auxout{\string
      \newlabel{#1}{{\@currentlabel}{\thepage}}}}}\@gtempa
   \if@nobreak \ifvmode\nobreak\fi\fi\fi
  \@esphack}

\def\alabel#1#2{\label{#1}\global\@fewtabfalse
    \ifcase\@eqcnt \def\@tempa{& & &}\or \def\@tempa{& &}
      \or \def\@tempa{&}
      \or\def\@tempa{}\fi\@tempa
{\hbox to 3cm{\phantom{\rm(\theequation,#2)}
\draftnote \hfil}\hskip -3cm {\rm(\theequation,#2)}}}

\def\clabel#1{\label{#1}\global\@fewtabfalse
    \ifcase\@eqcnt \def\@tempa{& & &}\or \def\@tempa{& &}
      \or \def\@tempa{&}
      \or\def\@tempa{}\fi\@tempa
{\hbox to 3cm{\phantom{\rm(\theequation)}
\draftnote \hfil}\hskip -3cm{\rm(\theequation)}}}

\def\eqnarray{\def\draftnote{{}}\global\@fewtabtrue
\stepcounter{equation}\let\@currentlabel=\theequation
\global\@eqnswtrue
\global\@eqcnt\z@\tabskip\@centering\let\\=\@eqncr
$$\halign to \displaywidth\bgroup\@eqnsel\hskip\@centering\@eqcnt\z@
  $\displaystyle\tabskip\z@{##}$&\global\@eqcnt\@ne
  \hskip 1\arraycolsep \hfil$\displaystyle{##}$\hfil
  &\global\@eqcnt\tw@ \hskip 1\arraycolsep
$\displaystyle\tabskip\z@{##}$
\hfil  \tabskip\@centering&\global\@eqcnt\thr@@\llap{##}\tabskip\z@
\cr}

\def\endeqnarray{\@@eqncr\egroup
      \global\advance\c@equation\m@ne$$\global\@ignoretrue}

\def\@eqnnum{\hbox to 3cm{\phantom{\rm(\theequation)} \draftnote
                         \hfil}\hskip -3cm {\rm(\theequation)}}

\def\@@eqncr{\let\@tempa\relax
    \ifcase\@eqcnt \def\@tempa{& & &}\or \def\@tempa{& &}
      \or \def\@tempa{&}
      \or\def\@tempa{}
\fi\@tempa
\if@eqnsw
\if@fewtab\@eqnnum\fi
\stepcounter{equation}\fi\global
\@eqnswtrue\global\@eqcnt\z@\global\@fewtabtrue\cr}

\def\draftcite#1{\ifnum\draftcontrol=1#1\else{}\fi}

\def\@lbibitem[#1]#2{\item{}\hskip -3cm \hbox to 2cm
{\hfil$\scriptstyle\draftcite{#2}$}\hskip
1cm[\@biblabel{#1}]\if@filesw
     {\def\protect##1{\string ##1\space}\immediate
      \write\@auxout{\string\bibcite{#2}{#1}}}\fi\ignorespaces}

\def\@bibitem#1{\item\hskip -3cm \hbox to 2cm
{\hfil $\scriptstyle\draftcite{#1}$}\hskip 1cm
\if@filesw \immediate\write\@auxout
       {\string\bibcite{#1}{\the\value{\@listctr}}}\fi\ignorespaces}

\def\nsection#1{\section{#1}\setcounter{equation}{0}}

\def\nappendixe{\def\thesection{A}\section*{Appendix }
\def\theequation{{A.\arabic{equation}}}
\def\theproposition{{A.\arabic{proposition}}}
\setcounter{equation}{0}
\setcounter{proposition}{0}}

\def\draftdate{\number\month/\number\day/\number\year\ \ \ \hourmin }

\global\def\draftcontrol{0}
\catcode`\@=12

\def\theequation{{\thesection.\arabic{equation}}}

\def\qq{\begin{eqnarray}}
\def\qqq{\end{eqnarray}}
\def\rx#1{~(\ref{#1})}
\def\ex#1{eq.\rx{#1}}
\def\eex#1{eqs.\rx{#1}}
\def\cx#1{~\cite{#1}}
\def\rw#1{~\ref{#1}}

\newlength{\shiftwidth}
\def\shift#1{&&\hbox to \shiftwidth{\hfill $\displaystyle#1$}}
\newlength{\sshiftwidth}
\def\sshift#1{\lefteqn{\hbox to
\sshiftwidth{\hfill$\displaystyle#1$}}}

\def\Im{\mathop{{\rm Im}}\nolimits}

\def\ie{{\it i.e.\ }}
\def\eg{{\it e.g.\ }}
\def\cf{{\it cf.\ }}
\def\rhs{{\it r.h.s.\ }}
\def\lhs{{\it l.h.s.\ }}

\def\Rhs{RHS\ }
\def\Rhsp{RHS}
\def\wrt{WRT\ }
\def\wrtp{WRT}

\def\ordH{{|H_1(M,\ZZ)|} }
\def\p{^{\prime}}
\def\pp{^{\prime\prime}}

\def\prosign{\mathop{{\rm sign}}\nolimits}
\def\sign#1{\prosign\left(#1\right)}

\def\spint#1{\int\limits^{+\infty}_{\scriptstyle -\infty \atop
[{#1}]}}
\def\intinf{ \int_{-\infty}^{+\infty} }

\def\promod{\mathop{{\rm mod}}\nolimits}
\def\mod#1{\;(\promod #1)}

\def\leg#1{ \{ #1 \}_K }

\def\pr#1#2{ \noindent{\em Proof of #1~\ref{#2}.} }

\def\zpk#1{ Z\p(#1;K) }
\def\zpmk{ \zpk{M} }
\def\zpxk{ \zpk{X} }
\def\zpnmk{ Z\p_{(N)}(M;K) }
\def\ztrnmk{ Z^{({\rm tr})}_{(N)}(M;K) }
\def\zmk{ Z(M;K) }
\def\zcmk{ Z\c(M;K) }

\def\ztrk#1{ Z^{({\rm tr})}(#1;K) }
\def\ztrmk{ \ztrk{M} }
\def\ztrxk{ \ztrk{X} }

\def\spq{ \check{q} }
\def\xp{ p }
\def\xq{ q }
\def\xt{ t }
\def\xh{ h }
\def\s{ \ast }
\def\ta{ \tilde{a} }
\def\l{ \lambda }

\def\c{ ^{(c)} }

\def\sqj{\sq_j}
\def\spj{\sp_j}
\def\sp{ \xp^\s }
\def\sq{ \xq^\s }
\def\qsn{ \sq_{(n)} }

\def\anmk{ a_n(M;K) }

\def\etpik{ e^{2\pi i\over K} }

\def\zq{ \ZZ[\spq] }
\def\zb#1{ \ZZ\left[ {#1} \right] }
\def\zti{ \zb{\xt,\xt^{-1}} }

\def\mnN{ m+n \leq N }

\def\snzkt{ \sum_{n=0}^{K-2} }
\def\snzi{ \sum_{n=0}^{\infty} }
\def\snzN{ \sum_{n=0}^N }
\def\sjN{ \sum_{j=1}^N }
\def\pjN{ \prod_{j=1}^N }
\def\pjtN{ \prod_{j=3}^N }
\def\sjtN{ \sum_{j=3}^N }
\def\smo{ \sum_{\mu_1=\pm 1} \mu_1 }
\def\smu{ \sum_{\mu=\pm 1} \mu }
\def\smz{ \sum_{m\geq 0} }
\def\snz{ \sum_{n\geq 0} }
\def\szmnt{ \sum_{0\leq m \leq {n\over 2} } }
\def\smps{ \sjN \mu_j \spj }
\def\slzm{ \sum_{l=0}^m }
\def\skztm{ \sum_{k=0}^{2m} }
\def\scomb{ \slzm \skztm }
\def\szns{ \sum_{0\leq n\leq \sjtN\spj-N+2} }
\def\smnN{ \sum_{\mnN} }
\def\sak{ \sum_{0\leq \a < K} }

\def\lcw{ \lambda_{\rm CW} }

\def\sg#1#2{ \sum_{0\leq #1 < 2K \atop #2 \in 2\ZZ + 1} }
\def\sgb{ \sg{\b}{\b} }
\def\sgbp{ \sg{\b\p}{\b} }
\def\sga{ \sg{\a}{\a} }
\def\sgap{ \sg{\a\p}{\a} }

\def\slnmhc{ \snzi \lnmc \xh^n }
\def\stlnmhc{ \snzi \tlnmc \xh^n }
\def\slnmh{ \snzi \lnm \xh^n }
\def\stlnmh{ \snzi \tlnm \xh^n }

\def\qlcwmc{ \spq^{ \left( 3\lcw(M) \right)^\vee} }
\def\qlcwm{ \spq^{ 3\lcw(M) } }
\def\qlcwx{ \spq^{ 3\lcw(X) } }
\def\tlcwc{ \left( 3\lcw \right)^\vee }

\def\pfq{ {\xp \over \xq} }
\def\pfqp{ \left( \pfq \right) }
\def\pfhp{ \left( {P\over H} \right) }

\def\lnm{ \lambda_n(M) }
\def\tlnm{ \tilde{\lambda}_n(M) }
\def\lnmc{ \lambda_n^\vee(M) }
\def\tlnmc{ \tilde{\lambda}_n^\vee(M) }

\def\tlnx{ \tilde{\lambda}_n(X) }
\def\lnn{ \lambda_n^{(N)} }
\def\tlnn{ \tilde{\lambda}_n^{(N)} }
\def\lnmn{ \lnn(M) }
\def\tlnmn{ \tlnn(M) }

\def\dappk{ \Delta_A\pp(\cK) }

\def\denm{ \Delta_n(M) }

\def\xla{ X\left( -{p\over r}, {q\over s}, l+pq \right) }

\def\tpifk{ {2\pi i\over K} }

\def\pfkp{ \left( {\pi\over K} \right) }
\def\stfk{ \sqrt{2\over K} }

\def\ordHx{{|H_1(X,\ZZ)|} }

\def\nfm{ {\ordH^{3\over 2} \over \stfk \sin \pfkp} }
\def\nfmx{ {\ordHx^{3\over 2} \over \stfk \sin \pfkp} }

\def\xpq{ X\left( {\xp_1 \over \xq_1}, \ldots,
     {\xp_N \over \xq_N} \right) }

\def\pfqj{ {\xp_j \over \xq_j} }
\def\qfpj{ {\xq_j \over \xp_j} }

\def\zjN{ 0 \leq j \leq N}
\def\ojN{ 1 \leq j \leq N}
\def\tjN{ 3 \leq j \leq N}

\def\shp{ \sign{H\over P} }
\def\sipq{ \sign{\xp \over \xq} }

\def\pfq{ {\xp \over \xq} }
\def\qfp{ {\xq \over \xp} }
\def\psq{ \sp \xq }
\def\qsp{ \sq \xp }

\def\qga{ \spq^{4^\s \qsp \a^2 } }
\def\qgap{ \spq^{4^\s \qsp {\a\p}^2 } }
\def\qgfa{ \spq^{ {1\over 4}\pfq \a^2} }
\def\qgfap{ \spq^{ {1\over 4}\pfq {\a\p}^2} }
\def\qgam{ \spq^{-\a^2} }

\def\qtsv{ \spq^{-3\psq\dappk - 2^\s\mu\sp + 2^\s} }
\def\qts{ \spq^{3s(\xp,\xq) - {1\over 4}\pfq } }
\def\qtsf{ \spq^{-3{\xq\over\xp}\dappk-{1\over 2}{\mu\over
     \xp}+{1\over 2} } }

\def\numqa{ \qvdp{2^\s\sq\a} }
\def\numf{ \qvdp{ {1\over 2} {\a\over \xq} } }

\def\ipf{ {i\pi\over 4} }
\def\iptf{ i\pi {3\over 4} }
\def\ipft{ {i\pi\over 2} }
\def\iptk{ {i\pi\over 2K} }
\def\ipkt{ {i\pi K\over 2} }

\def\eipfkshp{ e^{\ipf \kappa \shp} }
\def\eiptfkshp{ e^{\iptf \kt \shp} }
\def\eipmfkspq{ e^{-\ipf \kappa \sipq} }
\def\eiptfkspq{ e^{-\iptf \kt \sipq} }
\def\eiko{ e^{\ipf (\kappa + 1)\sipq} }
\def\eikoo{ e^{\ipf (\kappa + 1)\left(\sipq+1 \right)} }

\def\emipts{ e^{-\ipft \shp} }

\def\eipfk{ e^{\ipf(\kappa - 1)} }
\def\eipfok{ e^{\ipf(1 - \kappa)} }
\def\eipfks{ e^{\ipf(\kappa - 1)\shp} }

\def\kt{ {K-2\over K} }

\def\psh{ P^\s H}
\def\hsp{ H^\s P}

\def\fph{ 4^\s\psh - 3 \sjN s^\vee(\xq_j, \xp_j) }
\def\fphn{ 4^\s(N-1)(N-2)\hsp - 2^\s(N-2)\hsp\smps +
  2^\s\hsp\sum_{1\leq j < j\p \leq N} \mu_j\mu_{j\p}\spj\xp^\s_{j\p}
  + 2^\s }
\def\ffphn{ {P\over 2H} \left( {(N-1)(N-2) \over 2} +
     (N-2)\sjN\mpj + \sum_{1\leq j<j\p \leq N}\mpj
     {\mu_{j\p}\over \xp_{j\p} } \right) }
\def\ffphnl{ P {(N-1)(N-2) \over 2} +
     (N-2) P \sjN\mpj + P \sum_{1\leq j<j\p \leq N}\mpj
     {\mu_{j\p}\over \xp_{j\p} } }

\def\mpj{ {\mu_j\over \xp_j} }

\def\fss{ \left( {H\over P} - 3\shp - 12\sjN s(\xq_j,\xp_j) \right) }
\def\fsse{ \left[{H\over P} - 3\shp - 12\sjN s(\xq_j,\xp_j) +
{P\over H} \smpnp^2 \right] }

\def\clkm{ C_{l,k;m} }
\def\cn{ C_n }

\def\qgb{ \spq^{-4^\s\psh\b^2} }
\def\qgbp{ \spq^{-4^\s\psh{\b\p}^2} }
\def\qgfb{ \spq^{-{1\over 4} {H\over P} \b^2} }
\def\qgfbp{ \spq^{-{1\over 4} {H\over P} {\b\p}^2} }

\def\numoN#1{ \pjN \qvdp{#1} }
\def\denb#1{ \qvdp{#1}^{N-2} }
\def\numeN#1{ \pjN \evdp{#1} }
\def\deneN#1{ \evdp{#1}^{N-2} }
\def\evdp#1{ \left( e^{-#1} - e^{#1} \right) }
\def\evdmp#1{ \left( e^{\pi{#1}} - e^{-\pi{#1}} \right) }

\def\freb{ {\numeN{i\pi{b\over\xp_j}} \over \deneN{i\pi b} } }

\def\qvdp#1{ \left( \qvd{#1} \right) }
\def\qvd#1{ \spq^{-#1} - \spq^{#1} }

\def\tvdt#1{ \xt^{-{#1\over 2} } - \xt^{#1\over 2} }

\def\tmvda{ \xt^{\a\over 2} - \xt^{-{\a\over 2}} }
\def\tmvdap{ \left( \tmvda \right) }

\def\qtds{ {1\over \qvd{2^\s} } }
\def\qatm{ \left( \spq^{\a\over 2}-\spq^{-{\a\over 2}}\right)^{2m}}

\def\tse{ T_\s(\epsilon) }
\def\tseN{ T^{N-2}_\s(\epsilon) }
\def\ftse{ \frac{ \pjtN \left( 1 - (1-\tse)^{\spj} \right) }{\tseN} }
\def\te{ T(\epsilon) }
\def\fte{ \frac{ \pjtN \left( 1 - (1-\te)^{1\over \xp_j} \right) }
     {T^{N-2}(\epsilon) }   }
\def\ftsea#1#2{ {1 - \left( 1 - #1 \right)^{#2} \over #1 } }

\def\smpsnp{ \left( \smps - N + 2 \right) }
\def\smpnp{ \left( \sjN\mpj - N + 2 \right) }

\def\qof{ {1\over 1 - \spq } }

\def\derek{ {1\over k!} \partial_{\epsilon}^{(k)} }
\def\dereka#1{ \left. \derek #1 \right|_{\epsilon=0} }
\def\derbm#1{ \left. \partial^{(2m)}_b #1 \right|_{b=0} }

\def\frqs{ \frac{ \pjtN \left( 1 - \spq^{\spj\b} \right)}
     {\left(1 - \spq^\b \right)^{N-2}} }
\def\frq{ \frac{ \pjtN \left( 1 - \spq^{\b\over \xp_j} \right)}
     {\left(1 - \spq^\b \right)^{N-2}} }

\def\frbq{ {\numoN{ {1\over 2}{\b\over \xp_j} } \over
       \denb{ {\b\over 2} } }  }

\def\hspcm{ {\hsp n^2 \choose m} }
\def\hspl{ {\hsp \choose l} }

\def\pfhcl{ {{P\over H} \choose l} }

\def\shp{ \sign{H\over P} }

\def\fgenx{ \nfmx \ztrxk }
\def\fgenm{ \nfm \ztrmk }
\def\fgennm{ \nfm \ztrnmk }

\def\lxom{ \lim_{x\rightarrow 1^-} }

\def\snL{ \sum_{n\in \Lambda} }

\def\apfh{ \left|{P\over H}\right|^{1\over 2} }

\def\vp#1{ \left[ #1 \right]^\vee }

\def\jakq{ J_\a(\cK;\spq) }
\def\jakt{ J_\a(\cK;\xt) }
\def\vakq{ V_\a(\cK;\spq) }
\def\vakt{ V_\a(\cK;\xt) }
\def\vnakt{ V^{(N)}_\a(\cK;\xt) }
\def\vnakq{ V^{(N)}_\a(\cK;\spq) }
\def\junkt{ J_\a({\rm unknot};\xt) }
\def\junkq{ J_\a({\rm unknot};\spq) }

\def\apfqr{ \left| \pfq \right|^{1\over 2} }
\def\aqfpr{ \left| \qfp \right|^{1\over 2} }

\def\DmnK{ D_{m,n}(\cK) }
\def\dmnK{ \dmn(\cK) }
\def\dmn{ d_{m,n} }

\def\qspqmk{ \spq^{-\psq(m-k)^2 + \sp(1+\mu\xq)(m-k)} }

\def\kf#1{ \chi_#1(4_1) }
\def\ks#1{ \chi_#1(6_1) }

\begin{document}

\begin{titlepage}
\centerline{\hfill                 q-alg/9601015}
\vfill
\begin{center}

{\large \bf
On $p$-Adic Convergence of Perturbative Invariants of Some Rational
Homology Spheres.
}
\\

\bigskip
\centerline{L. Rozansky
}

\centerline{\em School of Mathematics, Institute for Advanced Study}
\centerline{\em Princeton, NJ 08540, U.S.A.}
\centerline{{\em E-mail address: rozansky@math.ias.edu}}

\vfill
{\bf Abstract}

\end{center}
\begin{quotation}
R.~Lawrence has conjectured that for rational homology spheres, the
series of Ohtsuki's invariants converges $p$-adicly to the $SO(3)$
Witten-Reshetikhin-Turaev invariant. We prove this
conjecture for Seifert rational homology spheres. We also derive it
for manifolds constructed by a surgery on a knot in $S^3$. Our
derivation is based on a conjecture about the colored Jones
polynomial that we have formulated in our previous paper. We also
present numerical examples of $p$-adic convergence for some simple
manifolds.
\end{quotation}
\vfill
\end{titlepage}

\pagebreak

\nsection{Introduction}
\label{s1}
\hyphenation{Re-she-ti-khin}
\hyphenation{Tu-ra-ev}

The cyclotomic properties of the
 $SO(3)$ Witten-Reshetikhin-Turaev
(\wrtp) invariant
\linebreak
$\zpmk$ of $3d$ manifolds $M$ have attracted a lot
of attention recently. This invariant was defined by R.~Kirby and
P.~Melvin\cx{KM} by modifying the Reshetikhin-Turaev surgery
formula of\cx{RT}. H.~Murakami showed that if $M$ is a rational
homology sphere (\Rhsp) and $K$ is an odd prime number, then
\qq
\zpmk \in \zq, \qquad \spq = \etpik,
\label{1.1}
\qqq
here $\zq$ is a cyclotomic ring:
\qq
{\spq^K - 1 \over \spq - 1} = 0
\label{1.01}
\qqq
(alternative proofs of\rx{1.1} were presented in\cx{MR} and \cx{Ro5}).

As an element of $\zq$, $\zpmk$ can be presented as a polynomial in
$\spq - 1$:
\qq
\zpmk = \snzkt \anmk \,\xh^n, \qquad
\xh = \spq - 1.
\label{1.2}
\qqq
The numbers $\anmk \in \ZZ$ depend on both the \Rhs $M$ and the
`level' $K$. T.~Ohtsuki showed in\cx{Oh1} and\cx{Oh2} how to
reprocess $\anmk$ into the invariants of $M$ which are independent of
$K$.

We introduce the following notations. Since our prime number is $K$,
we use the term `$K$-adic' instead of the usual term `$p$-adic'. For
$\xq \in \ZZ$, $\xq \neq 0 \mod{K}$ let $\sq$ denote $K$-adic
inverse of $\xq$. In other words, $\xq^\s$ is a formal series
in positive powers of $K$
\qq
\sq = \snzi \qsn K^n, \qquad
0 \leq \qsn \leq K-1,
\label{1.3}
\qqq
such that for any $N\geq 0$
\qq
\xq \snzN \qsn K^n = 1 \mod{ K^{N+1} }.
\label{1.04}
\qqq
We denote by $\vee$ an operation that converts rational numbers
whose denominators are not divisible by $K$, into $K$-adic numbers:
\qq
\pfqp^\vee = \xp\sq, \qquad
\xp, \xq \in \ZZ, \;\; \xq \neq 0 \mod{K}.
\label{1.03}
\qqq
%
\begin{theorem}[Ohtsuki\cx{Oh1},\cx{Oh2}]
\label{t1.1}
Let $M$ be a \Rhsp. Then there exists an infinite sequence of
rational invariants $\lnm$, $n\geq 0$ such that if $K$ is an odd
prime and $K > \ordH$, $\ordH$ being the order of first homology,
then
\qq
\leg{\ordH} \ordH
\anmk = \lnmc \mod{K}\;\;\;{\rm for}\;\;
n \leq {K-3\over 2},
\label{1.4}
\qqq
here $\leg{\ordH}$ is the Legendre symbol: for $p\in \ZZ$,
$\leg{p}=1$ if there exists some $q\in \ZZ$ such that
$p = q^2\mod{K}$, $\leg{p}=-1$ otherwise.
\end{theorem}

R.~Lawrence suggested that this theorem can be strengthened:
\begin{conjecture}[Lawrence\cx{Lw}]
\label{c1.1}
If $M$ is an integer homology sphere then $\lnm \in \ZZ$ and the
cyclotomic series $\slnmh$ converges $K$-adicly:
\qq
\slnmh = \zpmk.
\label{1.9}
\qqq
If $M$ is a rational homology sphere then
\qq
\lnm \in \zb{ {1\over 2}, {1\over \ordH} }
\label{1.10}
\qqq
and
if $\ordH \neq 0 \mod{K}$ then the cyclotomic series $\slnmhc$
converges $K$-adicly:
\qq
\slnmhc = \leg{\ordH} \ordH \zpmk.
\label{1.6}
\qqq
\end{conjecture}
R.~Lawrence proved her conjecture for the subclass
\qq
\xla, \qquad l,p,q,r,s \in \ZZ, \qquad q=2, \qquad ps - qr = 1
\label{1.11}
\qqq
of 3-fibered Seifert rational homology spheres. The
manifolds\rx{1.11} were constructed by Dehn's surgery on a $(p,q)$
torus knot with $q=2$ and framing $l$.

Let us comment briefly on the notion of $K$-adic limit in the
cyclotomic ring\footnote{I am thankful to J.~Bernstein for explaining
to me the relevance of the notion of cyclotomic $K$-adic convergence
for interpreting the results and conjectures of\cx{Lw}.}. A
cyclotomic relation\rx{1.01} tells us that
\qq
\xh^{K-1} = - K \snzkt {1\over K} {K \choose n+1} \xh^n.
\label{1.7}
\qqq
If K is prime then
${1\over K} {K \choose n+1} \in \ZZ$ for $ 0\leq n \leq K-2$.
Therefore \ex{1.7} allows us to reduce any cyclotomic polynomial of
$\xh$ to the `fundamental' powers $\xh^n$, $0\leq n\leq K-2$. We
denote this reduction by $\spadesuit$. Since the coefficients of
$\left[ \xh^{n(K-1) + m} \right]^\spadesuit$, $m,n \geq 0$ are
divisible by $K^n$, we see that the map $\spadesuit$ converts the
`analytic' smallness of high powers of $\xh$ into $K$-adic smallness
of the fundamental coefficients. As a result, \ex{1.6} means that for
any $N_1>0$ there exists $N_2$ such that for any $N>N_2$, the
coefficients of the polynomial
\qq
\left[ \leg{\ordH} \ordH
\zpmk - \snzN \lnmc \xh^n \right]^\spadesuit
\label{1.8}
\qqq
are all divisible by $K^{N_1}$.

A certain amount of information about the first Ohtsuki invariants is
already known. H.~Murakami showed\cx{Mu1},\cx{Mu2} that
$\lambda_0(M)=1$ and $\lambda_1(M)=3\lcw(M)$, here $\lcw(M)$ is the
Casson-Walker invariant.
X-S.~Lin and Z.~Wang proved\cx{LW} that for any integer
homology sphere $M$, $\lambda_2(M)\in 3\ZZ$.

The Casson-Walker invariant is known to satisfy the
following properties\cx{AM},\cx{W}:
\qq
\lcw(M) \in 2\ZZ \qquad &\mbox{if}&\qquad \ordH = 1,
\label{2.y3}
\\
6\ordH \lcw \in \ZZ \qquad &\mbox{if}& \qquad \ordH > 1.
\label{2.y4}
\qqq
We conjecture that the fraction $1\over 2$ in\rx{1.10} is due
completely to\rx{2.y4}, so that both parts of Conjecture\rw{c1.1} can
be united in one statement:
\begin{conjecture}
\label{c1.2}
For a \Rhs $M$ there exists a sequence of rational invariants $\tlnm$
such that
\qq
\tlnm \in \zb{1\over \ordH},
\label{1.y5}
\qqq
and if $\ordH \neq 0 \mod{K}$ then the cyclotomic series $\stlnmhc$
converges $K$-adicly:
\qq
\qlcwmc \stlnmhc = \leg{\ordH} \ordH \zpmk.
\label{1.y6}
\qqq
\end{conjecture}
Conjecture\rw{c1.1} follows from this conjecture because of
the properties\rx{2.y3} and\rx{2.y4} of the Casson-Walker invariant.

In our previous paper\cx{Ro5} we showed\footnote{
The proof of\cx{Ro5} relies on the Reshetikhin formula which is a
special integral representation of the colored Jones polynomial of a
link. A `path integral proof' of this formula was presented
in\cx{Ro2}. We will give a mathematically rigorous proof in\cx{Ro?}.
It is exactly the same proof as in\cx{Ro5} except that we use
Kontsevich's integral formula\cx{Ko1} for the $1/K$ expansion of the
Jones polynomial rather than a Bar-Natan style\cx{BN} perturbation
theory (see Appendix of\cx{Ro6} for the sketch of this argument).}
that the formal power series $\snzi \lnm \xh^n$ coincides (up to a
normalization) with the `perturbative invariant' $\ztrmk$ which is a
generating function of invariants $\denm$ (see\cx{Ro6} and references
therein for the definition of $\ztrmk$ and $\denm$). In other words,
$\lnm$ appear as coefficients in the re-expansion of the formal power
series $\snzi \denm K^{-n}$ in powers of
$\xh = \tpifk + \cO(K^{-2})$ rather than $K^{-1}$. More precisely, in
the notations of\cx{Ro6},
\qq
\nfm \ztrmk
\label{1.12}
=
{\pfkp \over \sin \pfkp} \snzi \denm
K^{-n} & = & \snzi \lnm \xh^n
\\
& = &
\qlcwm \stlnmh.
\nonumber
\qqq

The relation\rx{1.12} tells us that we can calculate Ohtsuki's
invariants $\lnm$, $\tlnm$
with the help of the same surgery formula that we
used for $\denm$. If $M$ is constructed by a surgery on a link
$\cL\subset S^3$, then $\lnm$ and $\tlnm$
can be expressed explicitly\cx{Ro5} in
terms of the derivatives of the colored Jones polynomial of $\cL$.
These expressions lead us in Section~\ref{s2} to the proof of
Conjecture~\ref{c1.2} for Seifert rational homology spheres. They
also help us in Section~\ref{s3} to derive Conjecture~\ref{c1.2} for
\Rhs constructed by a surgery on a knot $\cK\subset S^3$ from the
`weak conjecture' of\cx{Ro7}. In Section\rw{s4} we present numerical
examples of Conjecture\rw{c1.1} for integral homology spheres.
In Section\rw{s5} we briefly discuss
this conjecture as a relation between the path integral and
number theoretical approaches to Witten-Reshetikhin-Turaev invariant.

\nsection{Seifert Rational Homology Spheres}
\label{s2}
In this section we will prove Conjecture~\ref{c1.2} for Seifert
rational homology spheres. A Seifert \Rhs $\xpq$ can be constructed
by a surgery on a framed $(N+1)$-component link $\cL\subset S^3$.
This link consists of the `primary' unknot and $N$ `secondary'
unknots simply connected to the primary unknot in a Hopf link style.
The primary unknot has zero framing while the secondary unknots have
rational framings $\pfqj$, $\zjN$ (this means that their parallels go
$\xp_j$ times around the unknots and $\xq_j$ times along the
unknots).

We introduce a notation
\qq
P = \pjN \xp_j, \qquad H = P \sjN \qfpj.
\label{2.1}
\qqq
$H$ determines the order of the first homology of the Seifert
manifold $X$:
\qq
\ordHx = |H|.
\label{2.2}
\qqq

Let $K$ be an odd prime number. In what follows we require in view of
\ex{2.2} that
\qq
H \neq 0 \mod{K}.
\label{2.3}
\qqq
This condition implies together with \ex{2.1} that there could be at
most one $\xp_j$ divisible by $K$. Therefore we assume that
\qq
\xp_j \neq 0 \mod{K}, \qquad 2\leq j \leq N.
\qqq
In fact, we will assume
for simplicity that
$\xp_j, \xq_j \neq 0 \mod{K}$, $\ojN$. This will allow us to use the
numbers $\spj, \sqj$ freely
in the intermediate equations. We claim however
that the only essential assumption is\rx{2.3}.

The formula for the $SO(3)$ invariant of the Seifert \Rhs $X$ was
derived in\cx{Ro5} (eq.(4.38)):
\qq
\zpxk & = & - {i\over 2\sqrt{K} } \eipfkshp \eiptfkshp
\label{2.4}
\\
&&\qquad
\times
\leg{|P|} \sign{P} \spq^{\fph} \qtds
\nonumber\\
&&\qquad
\times
\sgb \qgb {\numoN{2^\s\spj\b} \over \denb{2^\s\b} },
\nonumber
\qqq
here $s(\xp,\xq)$ is a Dedekind sum
\qq
s(\xp,\xq) = {1\over 4\xq} \sum_{j=1}^{\xq-1}
\cot\left( {\pi j\over \xq}\right)
\cot\left( {\pi \xp j\over \xq}\right)
\label{2.5}
\qqq
and
\qq
\kappa =
\left\{
\begin{array}{cl}
1 &\mbox{if $K=1 \mod{4}$}\\
-1&\mbox{if $K=-1 \mod{4}$}.
\end{array}
\right.
\label{2.6}
\qqq
In\cx{Ro5} $\s$ meant only inverse $\promod{K}$: $\xp\sp = 1\mod{K}$,
while in this paper $\s$ means $K$-adic inverse. This does not make
any difference in the context of \ex{2.4} because all `asterisk'
numbers are in the exponent of $\spq$, so only their value
$\promod{K}$ is important.

To prove Conjecture~\ref{c1.2} for the Seifert \Rhs $X$ it is better
to use a slightly different formula for its $SO(3)$ invariant.
\begin{lemma}
\label{l2.1}
The $SO(3)$ invariant $\zpmk$ can be presented as a $K$-adicly
converging series:
\qq
\zpxk & = &
-\leg{|H|} \sign{H} \smo
\label{2.22}
\\
&&\!\!\!\qquad
\times
\spq^{\tlcwc + \fphn}
\nonumber\\
&&\!\!\!\qquad
\times
\qof \smz \scomb \clkm {\hsp \choose l} \dereka{\ftse},
\nonumber
\qqq
here
\qq
\mu_j =
\left\{
\begin{array}{cl}
\pm 1 &\mbox{for $j=1$}\\
1&\mbox{for $2\leq j\leq N$},
\end{array}
\right.
\label{2.23}
\qqq
and $\lcw$ is the Casson-Walker invariant of $X$ as calculated by
C.~Lescop\cx{L}:
\qq
\lcw =
{1\over 12}{P\over H} \left( 2 - N + \sjN {1\over \xp_j^2} \right)
- {1\over 4} \sign{H\over P} + {1\over 12}{H\over P} -
\sjN s(\xq_j, \xp_j).
\label{2.18}
\qqq
The coefficients $\clkm \in \ZZ$ are defined by \ex{2.14} and
\qq
\tse = 1 - (1+\epsilon) \spq^{-\hsp\smpsnp}.
\label{2.x24}
\qqq
\end{lemma}

\pr{Lemma}{l2.1}
We derive \ex{2.22} from \ex{2.4}. First of all, we want to remove
the factor ${1\over 2}$ from the \rhs of \ex{2.4}. All the factors in
the sum over $\b$ have definite parity. The range of summation could
also be made symmetric ($-K\leq \b \leq K$, $\b\in 2\ZZ+1$) due to
the periodicity of the summand under the shift
$\b\rightarrow \b+K$. Therefore we can place
$2\spq^{-2^\s\xp_1^\s\b}$
instead of $\qvd{2^\s\xp_1^\s\b}$ without changing the sum.

Next, we are going to transform the summand of \ex{2.4} by
`completing the square' in the exponent:
\qq
\spq^{-4^\s\psh\b^2} \spq^{\pm 2^\s\spj\b} =
\spq^{-4^\s\psh(\b\mp\hsp\spj)^2 + 4^\s\hsp{\spj}^2}.
\label{2.7}
\qqq
We introduce a new variable $\b\p$:
\qq
\b = \b\p - \hsp \smpsnp
\label{2.8}
\qqq
($\mu_j$ are defined by \ex{2.23}) and rewrite \ex{2.4} as
\qq
\zpxk & = & - {i\over 2\sqrt{K} } \eipfkshp \eiptfkshp
\label{2.9}
\\
&&\qquad
\times
\leg{|P|} \sign{P} \spq^{\fph + 4^\s\hsp\smpsnp^2}
\nonumber\\
&&\qquad
\times
\qtds
\sgbp \qgbp \frqs.
\nonumber
\qqq
We kept the variable $\b$ in the last factor of this equation meaning
that it is a function\rx{2.9} of $\b\p$.

To calculate the sum over $\b\p$ we present the last fraction factor
of the \rhs of \ex{2.9} as an explicit polynomial in $\spq^{\b}$:
\qq
\frqs = \szns \cn \spq^{n\b}, \qquad
\cn \in \ZZ.
\label{2.10}
\qqq
The coefficients $\cn$ depend, of course, on the numbers $\spj$ and
$N$. The sum over $\b\p$ for individual monomials $\spq^{n\b}$ can be
calculated by completing the square:
\qq
\sgbp \qgbp \spq^{n\b} & = &
\sgbp \spq^{-4^\s\psh \left(\b\p - 2\hsp n\right)^2}
\spq^{\hsp n^2} \spq^{-\hsp \smpsnp n}
\label{2.11}
\\
& = &
\sqrt{K} \eipfk \leg{\psh} \spq^{\hsp n^2} \spq^{-\hsp\smpsnp n}.
\nonumber
\qqq

Let us substitute $\spq = 1 + \xh$ in $\spq^{\hsp n^2}$ and present
the latter as a $K$-adicly convergent series
\qq
\spq^{\hsp n^2} = (1+\xh)^{\hsp n^2} = \smz  \hspcm
\xh^m,
\label{2.12}
\qqq
here by definition
\qq
{x\choose m} = {1\over m!} \prod_{l=0}^{m-1} (x-l),
\qquad
x\in \IR, \;\; m\in\ZZ, \;\;m\geq 0.
\label{2.13}
\qqq
The polynomial $P(x,n) = {xn^2 \choose m}$ takes integer values when
$x,n\in \ZZ$. Therefore it can be presented as an integer linear
combination of the product of elementary binomial polynomials of $x$
and $n$:
\qq
{xn^2 \choose m} = \scomb \clkm {x\choose l}{n\choose k},
\qquad \clkm \in \ZZ.
\label{2.14}
\qqq
After presenting a binomial polynomial ${n\choose k}$ as a derivative
\qq
{n\choose k} = \dereka{(1+\epsilon)^n},
\label{2.15}
\qqq
we finally re-express $\spq^{\hsp n^2}$ as
\qq
\spq^{\hsp n^2} = \smz \xh^m \scomb \clkm \hspl
\dereka{(1+\epsilon)^n}.
\label{2.16}
\qqq
Now we combine \eex{2.10},\rx{2.11} and\rx{2.16} in order to
calculate the sum over $\b\p$:
\qq
\hspace{-0.5in}
\sgbp \qgbp \spq^{n\b}& = &
\sqrt{K} \eipfk \leg{\psh}\!\!\! \szns \hspace{-0.4in}\cn
\spq^{-\hsp\smpsnp n}
\label{2.17}
\\
&&\qquad
\times
\smz \xh^m \scomb \clkm \hspl \dereka{(1+\epsilon)^n}
\nonumber
\\
& = &
\sqrt{K} \eipfk \leg{\psh} \smz \xh^m \scomb \clkm \hspl
\nonumber
\\
&&\qquad
\times
\dereka{\ftse},
\nonumber
\qqq
the function $\tse$ was defined in \ex{2.x24}.
We used \ex{2.10} `backwards' in order to get rid of the sum over
$n$.
If we substitute the
formula\rx{2.17}
in \ex{2.9} and use \ex{2.18} together with simple relations
\qq
& 4^\s = { 1- \kappa K \over 4} \mod{K},
\label{2.19}
\\
& \eipfk \leg{|P|} \leg{\psh} = \eipfks \leg{|H|},
\label{2.20}
\\
& i \sign{P} \emipts = \sign{H}
\label{2.21}
\qqq
in order to simplify the exponentials in front of the sum over
$\b\p$, then we obtain \ex{2.22}. $\Box$

Note that the requirement $\xp_1 \neq 0 \mod{K}$ looks artificial
indeed, because whenever $\xp_1^\s$ appears in \ex{2.22}, it is
always canceled by a factor of $P$.

Now we turn to the perturbative invariant $\ztrxk$ which generates
Ohtsuki's invariants $\tlnx$ via \ex{1.12}. The formula for $\ztrxk$
was derived in\cx{RoS1},\cx{Ro1}. It has obvious similarities with
\ex{2.4}:
\qq
\ztrxk & = & - {1\over 4K\sqrt{|P|} } e^{\iptf\shp} \sign{P}
\spq^{ {1\over 4}\fss }
\label{2.022}
\\
&&\qquad
\times
\spint{\b=0} d\b\,\qgfb \frbq.
\nonumber
\qqq
The symbol $\spint{\b=0}$ means that we take only the contribution of
the stationary phase point $\b=0$ to the integral. In other words, we
have to expand the factor
\qq
\frbq
\qqq
in Taylor series in $\b$ at $\b=0$ and then integrate the series
together with the gaussian exponential $\qgfb$ term by term in order
to produce the expansion\rx{1.12}. The result would be
the series\rx{2.x1} (\cf eq.(1.8)
of\cx{RoS1}).

A more effective formula for the perturbative invariant follows from
eqs.~(3.21), (3.23) of\cx{RoS2} specialized to the trivial
connection. It can be obtained from \ex{2.022} by rotating the
integration contour in the complex plane by $\pi\over 4$ and
rescaling the integration variable.

\begin{proposition}
\label{p2.1}
The generating function $\fgenx$ is an analytic function of $K$ in
the area where $iK{H\over P}$ is not a positive real number:
\qq
\fgenx & = &
- {i\over 4} {H\over \sin\pfkp}  e^{\iptk \fss}
\label{2.y1}\\
&&\qquad\times
\intinf dx\, e^{-\pi x^2} F\left(e^{-{i\pi\over4} \sign{H\over P} }
\sqrt{ {2\over K} {P\over H} } \, x \right),
\nonumber
\\
F(b) & = & \freb
\label{2.y2}
\qqq
In the limit of large $K$ the integral can be approximated by
expanding the pre-exponential factor $F$ in powers of $x^2/K$ and
intergrating each term of the expansion separately with the gaussian
factor. This procedure results in the asymptotic series in $1/K$:
%
\qq
\fgenx & = &
- {i\over 4} {H \over \sin\pfkp} e^{\iptk \fss}
\label{2.x1}
\\
&&
\;\;\;\;
\times
\sum_{m\geq 1}
{1\over m!} \pfhp^m {1\over (2\pi iK)^m} \derbm{\freb}.
\nonumber
\qqq
\end{proposition}

The expansion\rx{2.x1} contains too many fractions with `large'
denominators, so we are going to construct a slightly different
expansion which would look similar to \ex{2.22}.
\begin{lemma}
\label{l2.2}
The perturbative invariant $\fgenx$ which generates Ohtsuki's
invariants $\tlnx$, can be presented as an asymptotic series in $\xh$:
\qq
\lefteqn{
\fgenx = \qlcwx \snzi \tlnx \xh^n
}
\label{2.24}
\\
& = &
-\sign{H} \smo \spq^{3\lcw(X) + \ffphn + {1\over 2} }
\nonumber
\\
&&
\qquad
\times
\qof \smz \xh^m \scomb \clkm \pfhcl \dereka{\fte},
\nonumber
\qqq
here
\qq
\te = 1 - (1+\epsilon)\,\spq^{ -{P\over H}\smpnp}
\label{2.024}
\qqq
\end{lemma}

\pr{Lemma}{l2.2}
The proof is absolutely similar to that of Lemma~\ref{l2.1}. We place
$2\spq^{-{1\over 2}{\b\over \xp_1} }$ instead of
$\qvd{ {1\over 2}{\b\over \xp_1} }$ and complete the square by
shifting the integration variable:
\qq
\b = \b\p - {P\over H}\smpnp,
\label{2.25}
\qqq
so that \ex{2.022} becomes
\qq
\ztrxk & = &
-{1\over 2K\sqrt{|P|} }e^{\iptf\shp} \sign{P} \smo
\nonumber
\\
&&
\qquad
\times
\spq^{ {1\over 4} \fsse}
\nonumber
\\
&&
\qquad
\times
\spint{\b\p=0} d\b\p\, \qgfbp \frq.
\label{2.26}
\qqq
In the last factor of the \rhs of this equation $\b$ means the
function\rx{2.25} of $\b\p$.

We present the last factor of the \rhs of \ex{2.26} as a
`generalized' geometric series in powers of $\spq^\b$:
\qq
\frq & = &
\lxom \left(\snzi x^n\spq^{n\b}\right)^{N-2}
\pjtN\left( 1-x^{1\over \xp_j}\spq^{\b\over \xp_j} \right)
\label{2.27}
\\
& = & \lxom \snL \cn x^n \spq^{n\b},
\qquad
\cn \in \ZZ,
\nonumber
\qqq
here $\Lambda$ is a certain set of non-negative rational numbers and
$\cn$ are multiplicities with which the powers $n$ appear in the
expansion. The individual powers $\spq^{n\b}$ are easy to
integrate:
\qq
\intinf d\b\p\, \qgfbp \spq^{n\b} =
\sqrt{2K} \apfh e^{-\ipf\shp} \spq^{ {P\over H} n^2 }
\spq^{ -n{P\over H} \smpnp}.
\label{2.28}
\qqq
We transform the quadratic exponential further:
\qq
\spq^{{P\over H}n^2} =
(1+\xh)^{{P\over H}n^2} =
\smz \xh^m \scomb \clkm \pfhcl {n\choose k},
\label{2.028}
\qqq
\ex{2.15} defines a binomial coefficient ${x\choose l}$ for rational
$x$.

Combining \eex{2.27},\rx{2.28} and\rx{2.028} we obtain the
following formula for the integral over $\b\p$:
\qq
\lefteqn{
\spint{\b\p=0} d\b\p\, \qgfbp \frq = \sqrt{2K} \apfh
}
\label{2.29}
\\
&&\times
e^{-\ipf\shp} \snL \cn \spq^{-n{P\over H}\smpnp}
\smz \xh^m \scomb \clkm \pfhcl \dereka{(1+\epsilon)^n}
\nonumber
\\
& = &
\sqrt{2K} \apfh e^{-\ipf\shp} \smz \xh^m \scomb \clkm \pfhcl
\dereka{\fte}.
\nonumber
\qqq
The function $\te$ is defined by \ex{2.024}. We used \ex{2.27}
backwards to get rid of the sum over $n$ in the middle expression of
\ex{2.29}.
After substituting the
expression\rx{2.29}
into \ex{2.26} and using \ex{2.18} we arrive at \ex{2.24}.
$\Box$

Now we are going to prove Conjecture~\ref{c1.2} with the help of two
lemmas:
\begin{lemma}
\label{l2.3}
For $\xq,m \in \ZZ$,
\qq
{{\xq\over \xp} \choose m} \in \zb{1\over \xp}
\label{2.30}
\qqq
and hence
\qq
& {{\xq\over \xp} \choose m}^\vee = {\sp \xq \choose m}
\qquad\mbox{for}\;\; \xp\in\ZZ,\;\;\xp\neq 0\mod{K},
\label{2.030}
\\
& \spq^{\xq\over \xp} = \smz { {\xq\over \xp} \choose m} \xh^m
\in \zb{1\over \xp}[[\xh]].
\label{2.0031}
\qqq
\end{lemma}
\begin{lemma}
\label{l2.4}
\qq
\ffphn + {1\over 2} \in \zb{1\over H}.
\label{2.031}
\qqq
\end{lemma}

\pr{Lemma}{l2.4}
It is obvious that
\qq
\ffphn + {1\over 2} \in \zb{ {1\over H}, {1\over 2} },
\label{2.1031}
\qqq
so\rx{2.031} is true if $H$ is even.

We assume now that $H$ is odd. It is sufficient to prove that
\qq
\ffphnl + H \in 2\ZZ.
\label{2.2031}
\qqq
The parity of all three terms in this expression is determined by the
parity of $N$ and individual numbers $\xp_j$. Since $H$ is
odd, there can be at most one even number among $\xp_j$.
Therefore\rx{2.2031} can be checked directly for all possible
combinations. $\Box$

\pr{Conjecture}{c1.2}
First, we prove \ex{1.6}. We compare the formulas\rx{2.22}
and\rx{2.24} with the help of \ex{2.0031}. It is obvious that
\qq
\lefteqn{
\vp{ \spq^{ \ffphn + {1\over 2} } }
}
\hspace{1.7in}
\label{2.32}
\\
& = &
\spq^{ 2^\s\hsp
\left( 2^\s(N-1)(N-2) - (N-2)\sjN \mu_j \spj +
\sum_{1\leq j<j\p \leq N} \mu_j\mu_{j\p} \spj \xp^\s_{j\p} \right)
+ 2^\s}
\nonumber
\qqq
Also,
\qq
\vp{ \spq^{-{P\over H}\smpnp} }
=
\spq^{-\hsp \smpsnp},
\label{2.34}
\qqq
hence
\qq
\vp{ \te } = \tse.
\label{2.35}
\qqq
Then we see that for $\tjN$
\qq
\lefteqn{
\vp{
\ftsea{\te}{1\over \xp_j}
}
=
-\vp{
\smz (-1)^m { {1\over \xp_j} \choose m+1} T^m(\epsilon)
}
}
\hspace{2in}
\label{2.36}
\\
&&
=
\smz (-1)^m {\spj \choose m+1}T_\s^m(\epsilon)
= \ftsea{\tse}{\spj}.
\nonumber
\qqq
Finally,
\qq
{ {P\over H}\choose l}^\vee = {\hsp \choose l}.
\label{2.33}
\qqq
Combining \eex{2.32},\rx{2.36} and\rx{2.33} we arrive at \ex{1.y6} for
the Seifert \Rhs $X$.

Now we turn to condition\rx{1.y5}. The invariants $\tlnx$ appear as
coefficients in the expansion of the \rhs of \ex{2.24}
(multiplied by $\spq^{-3\lcw(X) }$)
in powers of
$\xh$. Lemma~\ref{l2.3} tells us that these coefficients are rational
and their denominators contain only the divisors of $H$ and $\xp_j$
for $\tjN$:
\qq
\tlnx \in \zb{ {1\over H}, {1\over \xp_3}, \ldots, {1\over \xp_N} }
\label{2.37}
\qqq
In our derivation of \ex{2.24} we made an arbitrary selection of two
factors $\qvd{{1\over 2}{\b\over\xp_j} }$, $j=1,2$ to be absorbed by
the gaussian exponential. We could select any other pair of numbers
$j$ and obtain a modified condition\rx{2.37}. The formula\rx{2.1} for
$H$ indicates that if a number is a divisor of at least one of any
$N-2$ numbers among all $\xp_j$, then that number divides $H$. This
means that the intersection of rings\rx{2.37} corresponding to all
possible choices of two missing numbers $\xp_j$, is exactly
$\zb{ {1\over H} }$. This proves\rx{1.y5}.$\Box$

\nsection{Surgery on a Knot}
\label{s3}
In this section we will study rational homology spheres constructed
by a rational surgery on a knot $\cK\subset S^3$. We will derive
Conjecture~\ref{c1.2} for these manifolds from Conjecture~1.1
of\cx{Ro7}.

Let $M$ be a \Rhs constructed by a rational $(\xp,\xq)$ surgery on a
knot $\cK\subset S^3$. Since
\qq
\ordH = |\xp|,
\label{3.1}
\qqq
we assume that
\qq
\xp \neq 0 \mod{K}.
\label{3.02}
\qqq
We also assume for simplicity that $\xq\neq 0 \mod{K}$, but we claim
that the latter assumption is not essential.

The formula for the $SO(3)$ invariant of $M$ was derived in\cx{Ro5}
along the lines of\cx{Je} from the Kirby-Melvin\cx{KM} modification
of the Reshetikhin-Turaev\cx{RT} surgery formula:
\qq
\zpmk & = &
- {i\over 2\sqrt{K} } \leg{|\xq|} \sign{\xq} \eipmfkspq \eiptfkspq
\label{3.2}
\\
&&
\qquad
\times
\spq^{3s^\vee(\xp,\xq) - 4^\s\qsp}
\sga \qga \numqa \jakq,
\nonumber
\qqq
here $\jakt$ is the `unframed' (\ie reduced to zero self-linking
number) colored Jones polynomial of the knot $\cK$ (we use the
complex variable $\xt$ in order to distinguish it from the cyclotomic
variable $\spq = \etpik$). The Jones polynomial is normalized in such
a way that it is multiplicative under a disconnected sum and
\qq
\junkt = {\tvdt{\a} \over \tvdt{1} }.
\label{3.3}
\qqq
The Casson-Walker invariant of $M$ is calculated by the surgery
formula of\cx{W}:
\qq
\lcw(M) = -s(\xq,\xp) + {\xq\over \xp}\dappk,
\label{3.y1}
\qqq
here $\dappk$ is the second derivative of the Alexander polynomial of
$\cK$ at $z=0$, $\dappk \in 2\ZZ$.

It is convenient to use the Jones polynomial in a slightly different
normalization:
\qq
\vakt = {\jakt \over \junkt}, \qquad
\vakt \in \zb{\xt,\xt^{-1} }.
\label{3.4}
\qqq
After substituting $\vakq$ instead of $\jakq$ in \ex{3.1} and
shifting the summation variable $\a$ in order to complete the square
in the exponent, we come to the `semi-analog' of
Lemma~\ref{l2.1}.
\begin{lemma}
\label{l3.1}
The $SO(3)$ invariant of $M$ is expressed in terms of $\vakq$ as a
sum
\qq
\zpmk & = &
- {i\over \sqrt{K} } \leg{|q|} \sign{q} \eiko\,
\qlcwmc
\smu \qtsv
\nonumber
\\
&&
\qquad
\times
\qof \sgap \qgap \vakq.
\label{3.5}
\qqq
In this equation $\a$ denotes a function of $\a\p$:
\qq
\a = \a\p + \sp(1+\mu\xq)
\label{3.6}
\qqq
\end{lemma}

\pr{Lemma}{l3.1}
First, we use the symmetries of the summand and summation range in
\ex{3.2} in order to substitute $2\spq^{-2^\s\sq\a}$ instead of
$\qvd{2^\s \sq \a}$. Since $\a$ is kept odd in the sum of \ex{3.2} we
can also change $\spq^{\a\over 2}$ into $-\spq^{2^\s\a}$ in the
formula for $\junkq$, because
\qq
2^\s = {1 - \kappa K \over 2}.
\label{3.7}
\qqq
As a result,
\qq
\jakq = {\qvd{2^\s\a} \over \qvd{2^\s} } \vakq.
\label{3.8}
\qqq

The change of variables\rx{3.6} allows us to absorb the factors
$\spq^{-2^\s\sq\a}$ and $\spq^{\pm 2^\s\a}$ into the gaussian
exponential $\qga$. Finally, the relation
\qq
s(\xp,\xq) + s(\xq,\xp) = {\xp^2+\xq^2+1\over 12\xp\xq} -
{1\over 4}\sign{\xp\xq}
\label{3.9}
\qqq
leads us to \ex{3.5}. $\Box$

Now we turn to the perturbative invariant $\ztrmk$ which is a
generating function of Ohtsuki's invariants $\lnm$ according to
\ex{1.12}. The formula for $\ztrmk$ was derived in\cx{Ro1} (see
also\cx{Ro6} and references therein):
\qq
\ztrmk & = &
i {\sin\pfkp \over 2K\sqrt{|\xq|} }
\sign{\xq} \eiptfkspq \qts
\label{3.10}
\\
&&
\qquad
\times
\spint{\a=0} d\a\, \qgfa \numf \jakq.
\nonumber
\qqq
Here the symbol $\spint{\a=0}$ means that we take only the
contribution of the stationary phase point $\a=0$ to the integral.
This contribution is calculated in the form of formal power series in
$K^{-1}$ (or in $\xh = \spq - 1$) of \ex{1.12} by expanding the
colored Jones polynomial in Taylor series in $K^{-1}$ and $\a$ and
then integrating monomials of $\a$ individually with the gaussian
factor $\qgfa$. Note that since we assume the polynomial $\jakq$ to
be expanded in powers of $\a$, it does make sense to evaluate it not
only for $\a\in \ZZ$ but also for $\a\in \IR$.

By using $\vakq$ instead of $\jakq$, substituting
\qq
\a = \a\p + {1+\mu\xq\over \xp}
\label{3.11}
\qqq
and using \ex{3.9} we obtain a modified expression for the
perturbative invariant.
\begin{lemma}
\label{l3.2}
The generating function $\fgenm$ can be presented as a stationary
phase integral:
\qq
\fgenm & = &
{i\over \sqrt{2K} } \apfqr \sign{\xq} \e^{-\iptf\sipq }
\spq^{3\lcw(M)}
\smu
\label{3.12}
\\
&&
\qquad
\times
\qtsf \qof \spint{\a\p=0} d\a\p\, \qgfap \vakq,
\nonumber
\qqq
$\a$ in the \rhs of this equation denotes the function\rx{3.11}.
\end{lemma}

The proof of this lemma is completely similar to the proof of
Lemma~\ref{l3.1} and we drop it.

A simple form of the coefficients of Taylor expansion of $\vakt$ in
powers of $\xt-1$ allows us to calculate the sum of \ex{3.5} and the
integral of \ex{3.12} explicitly. We will prove Conjecture~\ref{c1.2}
for $M$ by comparing these expressions.

For a fixed value of $\a$, $\vakt \in \zb{\xt, \xt^{-1} }$. Therefore
we can expand $\vakt$ in Taylor series in powers of $\xt - 1$ at
$\xt=1$. According to the Melvin-Morton conjecture\cx{MM} which was
proven by D.~Bar-Natan and S.~Garoufalidis\cx{BG}, the coefficients
of this expansion are polynomials in $\a$ of a limited degree:
\qq
\vakt = \snz \szmnt \DmnK \a^{2m} (\xt-1)^n,
\label{3.13}
\qqq
here $\DmnK \in \IQ$ are Vassiliev invariants of the knot $\cK$ of
degree $n$. We can rearrange \ex{3.13} as an expansion in powers of
$\a(\xt-1)$ and $\xt-1$:
\qq
\vakt = \snz \smz D_{m,n+2m}(\cK)\; [\a(\xt-1)]^{2m} (\xt-1)^n.
\label{3.14}
\qqq
We go one step further and use a variable
\qq
\tmvda = \a(\xt-1) + \cdots
\label{3.15}
\qqq
instead of $\a(\xt-1)$:
\qq
\vakt = \snz \smz \dmnK \tmvdap^{2m} (\xt-1)^n.
\label{3.16}
\qqq
We made the following `weak' conjecture in\cx{Ro7}:
\begin{conjecture}
\label{c3.1}
All coefficients $\dmnK$ in the expansion\rx{3.16} are integer:
\qq
\dmnK \in \ZZ.
\label{3.17}
\qqq
\end{conjecture}

Consider an `approximate' Jones polynomial
\qq
\vnakt = \smnN \dmnK \tmvdap^{2m} (\xt-1)^n.
\label{3.18}
\qqq
It is easy to see that for a fixed $\a$, the coefficients of its
first $N+1$ terms in Taylor expansion in powers of $\xt-1$ coincide
with those of the `exact' polynomial $\vakt$. The following simple
lemma is an analog of Lemma~2.3 of\cx{Mu1}.
\begin{lemma}
\label{l3.3}
If $f(\xt)\in \zti$ and
$\left.\partial_\xt^{(n)} f(\xt)\right|_{\xt=0} = 0$ for
$0 \leq n \leq N$, then
there exists $g(\xt) \in \zti$ such that
$f(\xt)=(\xt-1)^{N+1}g(\xt)$.
\end{lemma}
It follows from this lemma that there exists a function
$g_\a(\xt)$ such that $g_\a(\xt) \in \zti$ for $\a\in \ZZ$ and
\qq
\vakt = \vnakt + (\xt-1)^{N+1} g_\a(\xt).
\label{3.19}
\qqq
Denote by $\zpnmk$ the \lhs of \ex{3.5} if we substitute $\vnakt$
instead of $\vakt$ in its \rhs.
\begin{lemma}
\label{l3.4}
The `approximate invariant' $\zpnmk$ belongs to the cyclotomic ring
$\zq$
\qq
\zpnmk \in \zq
\label{3.20}
\qqq
and converges $K$-adicly to the `exact' invariant $\zpmk$:
\qq
\lim_{N\rightarrow \infty} \zpnmk = \zpmk.
\label{3.21}
\qqq
\end{lemma}

\pr{Lemma}{l3.4}
To prove\rx{3.20} we calculate explicitly the sum of \ex{3.5} for an
individual monomial
$\dmn\qatm (\spq - 1)^n$:
\qq
{1\over \sqrt{K} } \sgap \qgap \qatm \!\!\!& = &
{1\over \sqrt{K} } \sgap \qgap \skztm {2m\choose k} \spq^{\a (m-k)}
\label{3.22}
\\
& = &
{1\over \sqrt{K} } \skztm {2m \choose k} \qspqmk
\sgap\qgap
\nonumber
\\
& = &
\eipfok \leg{\qsp} \skztm {2m \choose k} \qspqmk.
\nonumber
\qqq
Since
\qq
-i\leg{|\xq|} \sign{\xq} e^{\ipf(\kappa+1)\sipq} \eipfok \leg{\qsp}
= \sign{\xp} \leg{|p|},
\label{3.23}
\qqq
we conclude that
\qq
\zpnmk & = & { \leg{|p|} \over 1-\spq} \sign{\xp}
\qlcwmc
\smu
\qtsv
\label{3.24}
\\
&&
\qquad
\times
\smnN \dmnK\,\xh^n \skztm {2m\choose k}\, \qspqmk.
\nonumber
\qqq
The divisibility by $\xh = 1-\spq$ in the \rhs of this equation is
easy to see. Indeed, the coefficient at $\xh^0$ in the expansion of
the summand in the sum over $\mu$ can be obtained by substituting
there $\spq=1$. The resulting expression is independent of $\mu$
except for the overall prefactor of $\mu$, hence it is eliminated by
the sum over $\mu$. Thus \ex{3.24} demonstrates\rx{3.20}.

To prove the $K$-adic limit\rx{3.21} we recall that
\qq
\sak \qgam = \sqrt{K} \eipfk
\label{3.25}
\qqq
and
\qq
\sak \qgam = \xh^{ {K-1\over 2} } w,
\qquad w,w^{-1} \in \zq,
\label{3.26}
\qqq
hence
\qq
{i\over \sqrt{K} } \eiko = \eikoo {w^{-1}\over \xh^{K-1\over 2} }
= {w\p \over \xh^{K-1\over 2} },
\qquad
w\p\in \zq.
\label{3.27}
\qqq
Therefore if we substitute the `error term'
$(\spq - 1)^{N+1} g_\a(\spq)$ of \ex{3.19} into the \rhs of \ex{3.5}
instead of $\vakq$, then we find that
\qq
\zpmk - \zpnmk = \xh^{N - {K-1\over 2} } w^{\prime\prime},
\qquad
w^{\prime\prime} \in \zq.
\label{3.28}
\qqq
This relation implies $K$-adic limit\rx{3.21}. $\Box$

There is an `asymptotic' counterpart of Lemma~\ref{l3.4}. Denote by
$\ztrnmk$ an `approximate perturbative invariant' which results from
substituting $\vnakq$ instead of $\vakq$ in the \rhs of \ex{3.12}.
\begin{lemma}
\label{l3.5}
The approximate perturbative invariant $\ztrnmk$ converges `formally'
to the exact invariant $\ztrmk$:
\qq
\lim_{N\rightarrow \infty} \ztrnmk = \ztrmk.
\label{3.29}
\qqq
\end{lemma}
The `formal limit'\rx{3.29} means that for any $N\p >0$ there exists
$N_0$ such that for any $N\geq N_0$ the first $N\p$ approximate
invariants $\lnmn$, $0\leq n < N\p$,
extracted from $\ztrnmk$ via the analog of
\ex{1.12}, coincide with the exact invariants $\lnm$.

\pr{Lemma}{l3.5}
An easy estimate
\qq
{1\over \sqrt{K} } \intinf d\a\p\, \qgfap \qatm = \cO(\xh^m)
\label{3.30}
\qqq
demonstrates that an individual term
$\dmn \qatm \xh^n$ of the sum\rx{3.18} contributes only to
$\lambda_{n\p}(M)$, $n\p\geq n+m-1$. Therefore the terms with $\dmn$,
$m+n\geq N+1$ do not contribute to $\lnm$, $n\leq N-1$. As a result,
\qq
\lnmn = \lnm
\qquad {\rm for}\;\; n \leq N-1.
\label{3.31}
\qqq
This proves the Lemma. $\Box$

Our final lemma is
\begin{lemma}
\label{l3.6}
Approximate invariants $\zpnmk$ and $\ztrnmk$ satisfy
Conjecture~\ref{c1.2}. Namely, the `invariants' $\tlnmn$ extracted
from $\ztrnmk$ via the analog of \ex{1.12} have only `simple'
denominators
\qq
\tlnmn \in \zb{1\over \xp},
\label{3.32}
\qqq
and there is a $K$-adic relation
\qq
\vp{\fgennm}
= \leg{\ordH} \ordH \zpmk.
\label{3.33}
\qqq
\end{lemma}

\pr{Lemma}{l3.6}
An integral of \ex{3.12} can be easily evaluated for an individual
term of the sum\rx{3.18}:
\qq
\intinf d\a\p\, \qgfap \qatm
& = &
\intinf d\a\p \qgfap \skztm {2m\choose k} \spq^{\a(m-k)}
\label{3.34}
\\
& = &
\sqrt{2K} \aqfpr e^{\ipf \sipq } \skztm {2m\choose k}
\spq^{ {1+\mu\xq\over \xp}(m-k)} \spq^{-\qfp(m-k)^2}
\nonumber
\qqq
Since
\qq
i\sign{q} e^{-\ipft\sipq} = \sign{\xp},
\label{3.35}
\qqq
we find that
\qq
\fgennm
\label{3.36}
&=&
{\sign{\xq} \over 1 - \spq}
\,\spq^{3\lcw(M)}
\smu \qtsf
\\
&&\qquad\times
\smnN \dmnK \xh^n \skztm {2m\choose k}
\spq^{ {1+\mu\xq\over \xp}(m-k)} \spq^{-\qfp(m-k)^2}.
\nonumber
\qqq

We analyze the coefficients in the expansion of this formula in
powers of $\xh$ with the help of Lemma~\ref{l2.3}. The factor
$\qof = - {1\over \xh}$ does not produce negative powers of $\xh$
because the term at $\xh^0$ in the summand of $\sum_{\mu=\pm 1}$ is
killed by the sum over $\mu$ the same way as it happened in
\ex{3.24}. All the exponents of $\spq$ belong explicitly to
$\zb{1\over \xp}$ except the term
\qq
-{\mu\over 2\xp} + {1\over 2} = {1\over \xp} \cdot {\xp-\mu\over
2}.
\label{3.37}
\qqq
If $\xp$ is even then
${1\over 2} \in \zb{1\over\xp}$. If $\xp$ is odd then
$\xp - \mu \in 2\ZZ$. Thus we proved\rx{3.32}.

Equation\rx{3.33} comes from comparing the expansions of the \rhs
of \ex{3.24} and\rx{3.36} in powers of $\xh$ and applying \ex{2.030}.
$\Box$

The Conjecture~\ref{c1.2} for $M$ is a direct consequence of
Lemmas~\ref{l3.4}, \ref{l3.5} and \ref{l3.6}.

\nsection{Numerical Examples}
\label{s4}

In her paper\cx{L}, R.~Lawrence has already presented some numerical
evidence in support of Conjecture~\ref{c1.1} for the manifolds
constructed by a $(1,1)$ surgery on torus knots. The proof of
Conjecture~\ref{c1.2} in Section\rw{s3} for manifolds constructed by
a rational surgery on a knot, was based on another conjecture and
therefore deserves a numerical check. In fact, we verified the
underlying Conjecture\rw{c3.1} in\cx{Ro7} numerically with certain
precision for some simple knots. This result translates through the
proofs of Section\rw{s3} into numerical verification of
Conjecture\rw{1.2}. Here we present the examples of Ohtsuki and
$SO(3)$ invariants of some simple integer homology spheres for mainly
illustrative purposes.

Let $\chi_{\xq}(\cK)$ denote an integer homology sphere constructed
by a $(1,\xq)$ surgery on a knot $\cK\subset S^3$. We calculated
$SO(3)$ invariants at $K=5$
for six manifolds $\chi_{1,2,3}(4_1)$ and
$\chi_{1,2,3}(6_1)$ ($4_1$ is known as `8-knot', see \eg\cx{BZ} for
the pictures of knots; $\kf{1}$ coincides with the manifold
$M_{3,-4}$ of\cx{Lw}). The invariants are presented as
polynomials\rx{1.2}. Their coefficients $a_n$ are listed in
Table\rw{t4.1}. The Ohtsuki invariants $\lnm$, $0\leq n\leq 11$ of
the same manifolds were calculated with the help of \eex{1.12}
and\rx{3.12}. They are collected in Table\rw{t4.2}. The
Table\rw{t4.3} contains the coefficients $\ta_n$ defined by the
relation
\qq
\sum_{n=0}^3 \ta_n(M) \xh^n
= \left[ \sum_{n=0}^{11} \lnm \xh^n
\right]^\spadesuit.
\label{4.1}
\qqq
It is easy to check that
\qq
a_n = \ta_n \mod{5^3}
\label{4.2}
\qqq
in full agreement with $K$-adic limit\rx{1.9} of Conjecture\rw{1.1}.

\begin{table}
\qq
\begin{array}{l|*{6} {c} }
\hline
  & \kf{1} & \kf{2} & \kf{3} & \ks{1} & \ks{2} & \ks{3} \\
\hline
a_0 & 1    &    6   &   1    &    -4  &    1   &   1
\\
a_1 & 4    &    8   &   2    &    -7  &    1   &  -1
\\
a_2 & 4    &    5   &   3    &    -4  &    1   &   0
\\
a_3 & 1    &    1   &   1    &    -1  &    0   &   0
\end{array}
\nonumber
\qqq
\caption{The coefficients $a_n(\chi_{1,2,3}(4_1))$ and
$a_n(\chi_{1,2,3}(6_1))$.
}
\label{t4.1}
\end{table}
\begin{table}
\qq
\def\y#1{ {\scriptscriptstyle #1} }
\begin{array}{l|*{9} {c} }
\hline
&\l_0&\l_1&\l_2&\l_3&\l_4&\l_5&\l_6&\l_7&\l_8\\
\hline
\kf{1}&\y{1}&\y{-6}&\y{69}&\y{-1064}&\y{20770}&\y{-492052}&\y{13724452}
&\y{-440706098}&\y{16015171303}\\
\kf{2}&\y{1}&\y{-12}&\y{270}&\y{-8284}&\y{324109}&\y{-15440692}&
\y{867600594}&\y{-56182136200}&\y{4119997542641}\\
\kf{3}&\y{1}&\y{-18}&\y{603}&\y{-27684}&\y{1624005}&\y{-116107654}&
\y{9795346273}&\y{-952637958170}&\y{104938749980019}\\
\ks{1}&\y{1}&\y{-12}&\y{246}&\y{-6916}&\y{248171}&\y{-10848488}&
\y{559466999}&\y{-33256127501}&\y{2238888918356}\\
\ks{2}&\y{1}&\y{-24}&\y{996}&\y{-57200}&\y{4207360}&\y{-377586960}&
\y{40010000129}&\y{-4889051681203}&\y{676832345117585}\\
\ks{3}&\y{1}&\y{-36}&\y{2250}&\y{-195060}&\y{21677895}&\y{-2940578892}
&\y{471069858664}&\y{-87035716226366}&\y{18219886814495290}
\end{array}
\nonumber
\qqq
\qq
\def\y#1{ {\scriptstyle #1} }
\begin{array}{l|*{3} {c} }
\hline
&\l_9&\l_{10}&\l_{11}\\
\hline
\kf{1}&\y{-649815778392}&\y{29121224693198}&
\y{-1428607184648931}\\
\kf{2}&\y{-337497038076594}&\y{30545298789501813}&
\y{-3026955693949081520}\\
\kf{3}&\y{-12914345285817636}&\y{1756098474092255228}&
\y{-261481531613674066565}\\
\ks{1}&\y{-168382708511446}&\y{13992172490383855}&
\y{-1273139103575342212}\\
\ks{2}&\y{-104697132728500664}&\y{17896885984795441440}&
\y{-3350210179233727412993}\\
\ks{3}&\y{-4262018348304709814}&\y{1101770773206971015914}&
\y{-311911822166530321178763}
\end{array}
\nonumber
\qqq
\caption{Ohtsuki invariants $\lambda_n(\chi_{1,2,3}(4_1))$ and
$\lambda_n(\chi_{1,2,3}(6_1))$.
}
\label{t4.2}
\end{table}
\begin{table}
\qq
\def\y#1{ {\scriptscriptstyle #1} }
\hspace{-7in}
\lefteqn{
\begin{array}{l|*{4} {c} }
\hline
  & \ta_0 & \ta_1 & \ta_2 &\ta_3\\
\hline
\kf{1}&\y{-1993867921922235374}&\y{-2540697622186346871}&
\y{-1996846249561765121}&\y{-401007308958445124}\\
\kf{2}&\y{-4192821765514758215994}&\y{-5339478176138672075492}&
\y{-4195910356932699813245}&\y{-840880791373767359499}\\
\kf{3}&\y{-361301312804709430245624}&\y{-460018993992612540674998}&
\y{-361478222004980994064372}&\y{-72392944327565893476624}\\
\ks{1}&\y{-1764664468572475265504}&\y{-2247385397501328211257}&
\y{-1766080693688788025254}&\y{-353995059668565997251}\\
\ks{2}&\y{-4624501555889040993306124}&\y{-5887570615994061037881249}&
\y{-4626301765208877840239624}&\y{-926250467309322277387125}\\
\ks{3}&\y{-429983196405880892355573999}&
\y{-547364469010535537295989626}&\y{-430093801055893893139442375}&
\y{-86079592745692299905577000}
\end{array}
}
\nonumber
\qqq
\caption{The coefficients $\ta_n(\chi_{1,2,3}(4_1))$ and
$a_n(\chi_{1,2,3}(6_1))$.
}
\label{t4.3}
\end{table}
%

\nsection{Discussion}
\label{s5}
E.~Witten\cx{Wi1} has originally defined the $SU(2)$ \wrt invariant
of a 3d manifold $M$ as a path integral over the $SU(2)$ connections
$A$ on $M$:
\qq
\zmk = \int {\cal D} A\, e^{ {i(K-2)\over 4\pi} \int_M
\langle A,dA + {1\over 3}[A,A] \rangle},
\label{5.1}
\qqq
here $\langle \, , \, \rangle$ is an appropriately normalized Killing
form on $su(2)$. The methods of quantum field theory allow us to
calculate this integral by the stationary phase approximation for
large values of $K$. The stationary points of the exponent (which is
proportional to the Chern-Simons invariant) are flat connections.
Therefore the invariant $\zmk$ is presented in the large $K$ limit as
a sum over the contributions $\zcmk$ coming from connected components
of the moduli space of flat $SU(2)$ connections on $M$ ($c$ indexes
these components):
\qq
\zmk = \sum_c \zcmk.
\label{5.2}
\qqq
The individual contributions are expressed generally as a product of
a `classical exponential', a certain power of $K$ and an
asymptotically convergent series:
\qq
\zcmk = e^{ {iK\over 4\pi} S\c_{\rm CS} }
K^{ \dim H\c_1 - \dim H\c_0\over 2}
\snzi \Delta\c_n(M) K^{-n}.
\label{5.3}
\qqq
In this formula $S\c_{\rm CS}$ is the Chern-Simons invariant of
connections of $c$, $H\c_{0,1}$ are the cohomologies of $d$ twisted
by a connection of $c$ and the coefficients $\Delta\c_n(M)$ are
called `$(n+1)$-loop corrections'.

The asymptotic form\rx{5.1},\rx{5.2} of the \wrt invariant was tested
numerically and analyticly by D.~Freed, R.~Gompf\cx{FG},
L.~Jeffrey\cx{Je} as well as in\cx{RoS1} for lens spaces and Seifert
manifolds.

If $M$ is a rational homology sphere, then the trivial connection is
a separate point in the moduli space of flat connections on $M$.
Therefore the trivial connection should produce a distinct
contribution to the sum\rx{5.2}. We checked\cx{RoS1} that the
`perturbative invariant'\rx{2.022} represents such a contribution for
Seifert rational homology spheres. We have reasons to believe (see
path integral arguments in Section~3 of\cx{Ro1}) that this is a
general case: perturbative invariant $\ztrmk$ which is defined
through the surgery formula without any reference to path integrals
(see\cx{Ro6} and references therein), is equal to the trivial
connection contribution to the `total' invariant $\zmk$.
Therefore Conjecture~\ref{c1.1} leads to the
following speculation:

\noindent
{\em Up to a normalization, the same power series $\snz \lnm \xh^n$
converges asymptotically to the trivial connection contribution into
the $SU(2)$ \wrt invariant and converges $K$-adicly to the $SO(3)$
\wrt invariant}.

\hyphenation{Ro-berts}
\hyphenation{Thur-ston}
\section*{Acknowledgements}
I am thankful to D.~Bar-Natan, J.~Bernstein, I.~Cherednik,
P.~Deligne, G.~Felder, D.~Freed, S.~Garoufalidis, J.~Handfield,
L.~Kauffman, M.~Kontsevich, R.~Lawrence, P.~Melvin, J.~Roberts,
D.~Thurston, V.~Turaev, A.~Vaintrob and E.Witten for many useful
discussions.

This work was supported by the National Science Foundation
under Grant DMS 9304580.

\end{document}

\end{document}

\end{document}


\begin{thebibliography}{99}
%
%
%
\bibitem{AM}S.~Akbulut, D.~J.~McCarthy, {\em Casson's Invariant for
Oriented Homology 3-Spheres - an Exposition}, Mathematical Notes,
{\bf 36}, Princeton University Press, Princeton, 1990.
%
\bimn{BN}{D.~Bar-Natan}{Perturbative Chern-Simons Theory}{Jour. of
Knot Theory and its Ramifications}{4-4}{1995}{503-548}
%
\bibitem{BG} D.~Bar-Natan, S.~Garoufalidis, {\em On the
Melvin-Morton-Rozansky Conjecture}, preprint 1994.
%
\bibitem{BZ}G.~Burde, H.~Zieschang, {\em Knots}, de Gruyter, Berlin
and New York, 1985.
%
\bimn{FG}{D.~Freed, R.~Gompf}{Computer Calculation of Witten's
3-Manifold Invariant}{Commun. Math. Phys.}{141}{1991}{79-117}
%
\bimn{Je}{L.~Jeffrey}{Chern-Simons-Witten Invariants of Lens Spaces
and Torus Bundles, and the Semiclassical Approximation}{Commun. Math.
Phys.}{147}{1992}{563-604}
%
\bimn{KM}{R.~Kirby, P.~Melvin}{The 3-Manifold Invariants of Witten
and Reshetikhin-Turaev for $sl(2,C)$}{Invent.
Math.}{105}{1991}{473-545}
%
\bimn{Ko1}{M.~Kontsevich}{Vassiliev's Knot Invariants}{Adv. in Sov.
Math.}{16(2)}{1993}{137-150}
%
\bimn{Lw}{R.~Lawrence}{Asymptotic Expansions of
Witten-Reshetikhin-Turaev Invariants for Some Simple
3-Manifolds}{J.~Mod.~Phys.}{36}{1995}{6106-6129}
%
\bimn{L}{C.~Lescop}{Invariant de Casson-Walker des sph\`{e}res
d'homologie rationnelle fibr\'{e}es de Seifert}
{C.~R.~Acad.~Sci.~Paris}{310}{1990}{727-730}
%
\bibitem{LW} X-S.~Lin, Z.~Wang, {\em On Ohtsuki's Invariants of
Integral Homology Spheres, I}, preprint, q-alg/9509009.
%
\bibitem{MR}G.~Masbaum, J.~Roberts, {\em A Simple Proof of
Integrality of Quantum Invariants at Prime Roots of Unity}, preprint,
September 28, 1995.
%
\bimn{MM}{P.~Melvin, H.~Morton}{The Coloured Jones Function}{Commun.
Math. Phys.}{169}{1995}{501-520}
%
\bimn{Mu1}{H.~Murakami}{Quantum $SU(2)$-Invariants Dominate Casson's
$SU(2)$-Invariant}{Math. Proc. Camb. Phil. Soc.}{115}{1993}{253-281}
%
\bimn{Mu2}{H.~Murakami}{Quantum $SO(3)$-Invariants Dominate
the SU(2)-Invariant of Casson and Walker}{Math. Proc. Camb. Phil.
Soc.}{117}{1995}{237-249}
%
\bimn{Oh1}{T.~Ohtsuki}{A Polynomial Invariant of Integral Homology
3-Spheres}{Math. Proc. Camb. Phil. Soc.}{117}{1995}{83-112}
%
\bibitem{Oh2}T.~Ohtsuki, {\em A Polynomial Invariant of Rational
Homology 3-Spheres}, preprint UTMS 94-49, August 12, 1994.
%
\bimn{RT}{N. Reshetikhin, V. Turaev}{Invariants of 3-Manifolds via
Link Polynomials and Quantum Groups}
{Invent. Math.}{103}{1991}{547-597}
%
\bimn{RoS1}{L.~Rozansky}{A Large $k$ Asymptotics of Witten's
Invariant of Seifert
Manifolds}{Commun.~Math.~Phys.}{35}{1994}{5219-5246}
%
\bibitem{RoS2}L.~Rozansky, {\em Residue Formulas for the Large $k$
Asymptotics of Witten's Invariants of Seifert Manifolds. The Case of
$SU(2)$.}, preprint UMTG-179, hep-th/9412075, to appear in Commun.
Math. Phys.
%
%
\bimn{Ro1}{L.~Rozansky}{A Contribution of the Trivial Connection to
the Jones Polynomial and Witten's Invariant of 3d Manifolds
I}{Commun. Math. Phys.}{175}{1996}{275-296}
%
\bimn{Ro2}{L.~Rozansky}{Reshetikhin's Formula for the Jones
Polynomial of a Link: Feynman Diagrams and Milnor's Linking Numbers}
{Journ. Math. Phys.}{35}{1994}{5219-5246}
%
\bibitem{Ro5}L.~Rozansky, {\em Witten's Invariants of Rational
Homology Spheres at Prime Values of $K$ and Trivial Connection
Contribution}, preprint UMTG-183, q-alg/9504015, to appear in Commun.
Math. Phys.
%
\bibitem{Ro6}L.~Rozansky, {\em On Finite Type Invariants of Links and
Rational Homology Spheres Derived from the Jones Polynomial and
Witten-Reshetikhin-Turaev Invariant}, q-alg/9511025.
%
\bibitem{Ro7}L.~Rozansky, {\em Higher Order Terms in the
Melvin-Morton Expansion of the Colored Jones polynomial}, preprint,
q-alg/9601009.
%
\bibitem{Ro?}L.~Rozansky, in preparation.
%
\bimn{Wi1}{E.~Witten}{Quantum Field Theory and the Jones
Polynomial}{Commun. Math. Phys.}{121}{1989}{351-399}
%
\bibitem{W}K.~Walker, {\em An Extension of Casson's Invariant},
Annals of Mathematics Studies, {\bf 126}, Princeton University Press,
Princeton, 1992.
%
\end{thebibliography}

\begin{thebibliography}{99}
%
\bibitem{AL}M.~Alvarez, J.~M.~F.~Labastida, {\em Vassiliev Invariants
for Torus Knots}, preprint q-alg/9506009.
%
\bibitem{BN}D.~Bar-Natan, {\em On the Vassiliev Knot Invariants\/},
Topology, to appear.
%
\bibitem{BG} D.~Bar-Natan, S.~Garoufalidis, {\em On the
Melvin-Morton-Rozansky Conjecture}, preprint 1994.
%
\bimn{BL}{J.~S.~Birman, X-S.~Lin}{Knot polynomials and Vassiliev's
invariants}{Invent. Math.}{111}{1993}{225-270}
%
\bibitem{BZ}G.~Burde, H.~Zieschang, {\em Knots}, de Gruyter, Berlin
and New York, 1985.
%
\bimn{ILR}{J.~M.~Isidro, J.~M.~F.~Labastida,
A.~V.~Ramallo}{Polynomials for Torus Links from Chern-Simons Gauge
Theories}{Nucl. Phys.}{B398}{1993}{187-236}
%
\bimn{Je}{L.~Jeffrey}{Chern-Simons-Witten Invariants of Lens Spaces
and Torus Bundles, and the Semiclassical Approximation}{Commun. Math.
Phys.}{147}{1992}{563-604}
%
\bibitem{KL}L.~Kauffman, S.~Lins, {\em Temperley-Lieb Recoupling
Theory and Invariants of 3-Manifolds}
%
\bimn{MM}{P.~Melvin, H.~Morton}{The Coloured Jones Function}{Commun.
Math. Phys.}{169}{1995}{501-520}
%
\bimn{M}{H.~Morton}{The Colored Jones Function and Alexander
Polynomial for Torus Knots}{Math. Proc. Cam. Phil.
Soc.}{117}{1995}{129-135}
%
\bibitem{Oh2}T.~Ohtsuki, {\em A Polynomial Invariant of Rational
Homology 3-Spheres}, preprint UTMS 94-49, August 12, 1994.
%
\bibitem{Ro1}L. Rozansky, {\em A Contribution of the Trivial
Connection to the Jones Polynomial and Witten's Invariant of 3d
Manifolds I.}, preprint UMTG-172-93, UTTG-30-93, hep-th/9401061, to
appear in Commun. Math. Phys.
%
%
\bibitem{RoS2}L.~Rozansky, {\em Residue Formulas for the Large $k$
Asymptotics of Witten's Invariants of Seifert Manifolds. The Case of
$SU(2)$.}, preprint UMTG-179, hep-th/9412075.
%
%
\bibitem{Ro5}L.~Rozansky, {\em Witten's Invariants of Rational
Homology Spheres at Prime Values of $K$ and Trivial Connection
Contribution}, preprint UMTG-183, q-alg/9504015, to appear in Commun.
Math. Phys.
%
\bibitem{Ro6}L.~Rozansky, {\em On Finite Type Invariants of Links and
Rational Homology Spheres Derived from the Jones Polynomial and
Witten-Reshetikhin-Turaev Invariant}, q-alg/9511025.
%
\bibitem{Ro8}L.~Rozansky, {\em On $p$-adic Convergence of Perturbative
Invariants of some Rational Homology Spheres}, in preparation.
%
\end{thebibliography}

\begin{thebibliography}{99}
%
\bibitem{AS1} S. Axelrod, I. Singer, {\em Chern-Simons Perturbation
Theory}, Proceedings of XXth Conference on Differential Geometric
Methods in Physics (New York, 1991) (S.~Catto and A.~Rocha, eds)
World Scientific, 1992, 3-45.
%
\bimn{AS2}{S.~Axelrod, I.~Singer}{Chern-Simons Perturbation Theory II}
{Jour. Diff. Geom.}{39}{1994}{173-213}
%
\bibitem{BN1}D.~Bar-Natan, {\em Perturbative Aspects of the
Chern-Simons Topological Quantum Field Theory}, Ph.D. thesis,
Princeton University, 1991.
%
\bibitem{BN2}D.~Bar-Natan, {\em On the Vassiliev Knot Invariants\/},
Topology, to appear.
%
\bibitem{BN3} D. Bar-Natan, {\em Vassiliev Homotopy String Link
Invariants}, Harvard University preprint, 1993.
%
\bibitem{BG} D.~Bar-Natan, S.~Garoufalidis, {\em On the
Melvin-Morton-Rozansky Conjecture}, preprint 1994.
%
\bimn{BL}{J.~S.~Birman, X-S. Lin}{Knot polynomials and Vassiliev's
invariants}{Invent. Math.}{111}{1993}{225-270}
%
\bibitem{Ga} S.~Garoufalidis, {\em Relations among 3-Manifold
Invariants}, Ph.D. Thesis, University of Chicago, August 1992.
%
\bibitem{GO} S.~Garoufalidis, T.~Ohtsuki, {\em On Finite Type
Invariants V: Rational Homology 3-Spheres}, preprint.
%
\bimn{Je}{L.~Jeffrey}{Chern-Simons-Witten Invariants of Lens Spaces
and Torus Bundles, and the Semiclassical Approximation}{Commun. Math.
Phys.}{147}{1992}{563-604}
%
\bimn{Ko1}{M.~Kontsevich}{Vassiliev's Knot Invariants}{Adv. in Sov.
Math.}{16(2)}{1993}{137-150}
%
\bibitem{Li} X-S. Lin, {\em Milnor Link Invariants Are All of Finite
Type}, Columbia University preprint, 1992.
%
\bibitem{LW} X-S.~Lin, Z.~Wang, {\em On Ohtsuki's Invariants of
Integral Homology 3-Spheres, I\/}, preprint, q-alg/9509009.
%
\bimn{KM}{R.~Kirby, P.~Melvin}{The 3-Manifold Invariants of Witten
and Reshetikhin-Turaev for $sl(2,C)$}{Invent.
Math.}{105}{1991}{473-545}
%
\bibitem{Lw}R.~Lawrence, {\em Asymptotic Expansions of
Witten-Reshetikhin-Turaev Invariants for Some Simple 3-Manifolds\/},
preprint IHES/M/95/39, to appear in J.~Mod.~Phys.
%
\bibitem{MR}G.~Masbaum, J.~Roberts, {\em A Simple Proof of
Integrality of Quantum Invariants at Prime Roots of Unity}, preprint,
September 28, 1995.
%
\bibitem{MeMo} P.~Melvin, H.~Morton, {\em The Coloured Jones
Function\/}, Commun. Math. Phys., to appear.
%
\bimn{Mu1}{H.~Murakami}{Quantum $SU(2)$-Invariants Dominate Casson's
$SU(2)$-Invariant}{Math. Proc. Camb. Phil. Soc.}{115}{1993}{253-281}
%
\bimn{Mu2}{H.~Murakami}{Quantum $SO(3)$-Invariants Dominate
the SU(2)-Invariant of Casson and Walker}{Math. Proc. Camb. Phil.
Soc.}{117}{1995}{237-249}
%
\bimn{Oh1}{T.~Ohtsuki}{A Polynomial Invariant of Integral Homology
3-Spheres}{Math. Proc. Camb. Phil. Soc.}{117}{1995}{83-112}
%
\bibitem{Oh2}T.~Ohtsuki, {\em A Polynomial Invariant of Rational
Homology 3-Spheres}, preprint UTMS 94-49, August 12, 1994.
%
\bibitem{Oh3}T. Ohtsuki, {\em Finite Type Invariants of Integral
Homology 3-Spheres}, preprint UTMS 94-42, June 1, 1994.
%
\bibitem{Re}N.~Reshetikhin, private communications.
%
\bimn{RT}{N. Reshetikhin, V. Turaev}{Invariants of 3-Manifolds via
Link Polynomials and Quantum Groups}
{Invent. Math.}{103}{1991}{547-597}
%
\bibitem{Ro1}L. Rozansky, {\em A Contribution of the Trivial
Connection to the Jones Polynomial and Witten's Invariant of 3d
Manifolds I.}, preprint UMTG-172-93, UTTG-30-93, hep-th/9401061, to
be published in Commun. Math. Phys.
%
\bimn{Ro2}{L.~Rozansky}{Reshetikhin's Formula for the Jones
Polynomial of a Link: Feynman Diagrams and Milnor's Linking Numbers}
{Journ. Math. Phys.}{35}{1994}{5219-5246}
%
\bibitem{Ro3}L. Rozansky, {\em A Contribution of the Trivial
Connection to the Jones Polynomial and Witten's Invariant of 3d
Manifolds II.}, preprint UMTG-187-94, hep-th/9403021, to
be published in Commun. Math. Phys.
%
\bibitem{Ro4}L.~Rozansky,{\em The Trivial Connection Contribution to
Witten's Invariant and Finite Type Invariants of Rational Homology
Spheres}, preprint UMTG-182, q-alg/9503011, submitted to Commun.
Math. Phys.
%
\bibitem{Ro5}L.~Rozansky, {\em Witten's Invariants of Rational
Homology Spheres at Prime Values of $K$ and Trivial Connection
Contribution}, preprint UMTG-183, q-alg/9504015, submitted to Commun.
Math. Phys.
%
\bimn{Wi1}{E.~Witten}{Quantum Field Theory and the Jones Polynomial}
{Commun. Math. Phys.}{121}{1989}{351}
%
\end{thebibliography}
\end{document}